\documentclass[preprint]{elsarticle}

\usepackage{geometry} 
\usepackage{amsmath,amssymb,amsthm,graphicx,verbatim,bm,fancyhdr,bbm,enumerate}
\DeclareMathAlphabet{\mathpzc}{OT1}{pzc}{m}{it}  

\theoremstyle{definition}
\theoremstyle{theorem}
\newtheorem{theorem}{Theorem}[section]

\newtheorem{lemma}[theorem]{Lemma}
\newtheorem{corollary}[theorem]{Corollary}

\newtheorem{conjecture}[theorem]{Conjecture}
\newtheorem{remark}{Remark}
\newtheorem{definition}{Definition}
\theoremstyle{remark}

\renewcommand{\P}{\mathbb{P}}

\newcommand{\E}{\mathbb{E}}

\newcommand{\R}{\mathbb{R}}
\newcommand{\Z}{\mathbb{Z}}

\newcommand{\T}{\mathbb{T}}
\newcommand{\eps}{\varepsilon}

\newcommand{\set}[1]{\left\{#1\right\}} 
\newcommand{\ob}[1]{\left(#1\right)} 
\newcommand{\cb}[1]{\left[#1\right]} 
\newcommand{\abs}[1]{\left\vert#1\right\vert} 
\newcommand{\norm}[1]{\|#1\|}
\newcommand{\dnorm}[1]{\left | \left | \left |#1\right | \right |\right |}

\newcommand{\ip}[1]{\left \langle#1\right \rangle}

\begin{document}
\begin{frontmatter}

\title{On the well-posedness of the stochastic Allen-Cahn equation in two dimensions}

\author[rvt]{Marc D. Ryser\corref{cor1}\fnref{fn1}\fnref{fn2}}
\ead{ryser@math.mcgill.ca}
\author[focal]{Nilima Nigam}
\ead{nigam@math.sfu.ca}
\author[focal]{Paul F. Tupper}
\ead{pft3@math.sfu.ca}

\cortext[cor1]{Corresponding author. }
\fntext[fn2]{Phone: +1-604-312-7973, Fax: +1-778-782-4947}
\fntext[fn1]{Research also conducted at Department of Mathematics, Simon Fraser University, 8888 University Drive, Burnaby, British Columbia  V5A 1S6, Canada }

\address[rvt]{Department of Mathematics and Statistics, McGill University, 805 Sherbrooke Street West, Montreal, Quebec H3A 2K6, Canada}
\address[focal]{Department of Mathematics, Simon Fraser University, 8888 University Drive,
Burnaby, British Columbia  V5A 1S6, Canada}

\begin{abstract}
White noise-driven nonlinear stochastic partial differential equations (SPDEs) of parabolic type are frequently used to model physical systems in space dimensions $d=1,2,3$. Whereas existence and uniqueness of weak solutions to these equations are well established in one dimension, the situation is different for $d\geq 2$. Despite their popularity in the applied sciences, higher dimensional versions of these SPDE models are generally assumed to be ill-posed by the mathematics community. We study this discrepancy on the specific example of the two dimensional Allen-Cahn equation driven by additive white noise. Since it is unclear how to define the notion of a weak solution to this equation, we regularize the noise and introduce a family of approximations. Based on  heuristic arguments and numerical experiments, we conjecture that these approximations exhibit divergent behavior in the continuum limit.
The results strongly suggest that shrinking the mesh size in simulations of the two-dimensional white noise-driven Allen-Cahn equation does not lead to the recovery of a physically meaningful limit.

\end{abstract}


\end{frontmatter}
\vspace{1cm}
\section{Introduction}\label{intro}
Stochastic equations of the type
\begin{align}\label{e0}
\partial_t u= A u + f(u) + \xi,
\end{align}
where $A$ is a linear elliptic differential operator, $f$ is a nonlinear function and $\xi$ is space-time white noise, play a central role in the modeling of a whole variety of phenomena in the physical sciences. Prominent examples are the Swift-Hohenberg equation in the study of thermal convection \cite{swift1977hydrodynamic} or the Kardar-Parisi-Zhang (KPZ) and Lai-Das Sarma-Villain equations in surface growth \cite{kardar1986dynamic,lai1991kinetic}. Another important domain of application for nonlinear parabolic SPDEs driven by additive white noise is the theory of the dynamics of critical phenomena \cite{goldenfeld1992lectures,chaikin2000principles,hohenberg1977theory}. Hohenberg and Halperin \cite{hohenberg1977theory} introduced a nowadays widely used classification of the various models of dynamic critical phenomena. Their classification includes among others: {\it Model A} for non-conserved quantities, e.g.\ the stochastic time-dependent Ginzburg-Landau and Allen-Cahn equations; {\it Model B} for conserved quantities, e.g.\ the stochastic Cahn-Hilliard equation; {\it Model C} which couples conserved and non-conserved fields, e.g.\ phase-field models of eutectic growth.  \\

\noindent Over the past decade, the applied science community has  paid a lot of attention to the two- and three-dimensional versions of these white noise-driven models. In addition to some analytic work, see e.g.\ \cite{wio2011recent} for the KPZ equation, emphasis has been put on numerical investigations: standard finite difference and pseudospectral methods for spatial discretization have been combined with Euler-Maryama or stochastic Runge-Kutta schemes to find numerical approximations to these SPDEs. For specific examples of such work we refer to \cite{ma2007scaling,giada2002pseudospectral} for the KPZ equation, to \cite{benzi2011phase,slutsker2008phase,ibanes2000dynamics,oguz1990domain,rao1993kinetics} for {\it Models A} and {\it B}, and to  \cite{ni2009chessboard,granasy2004general,granasy2002nucleation,elder1994stochastic,drolet2000phase} for {\it Model C}. \\

\noindent But while these models are extensively studied by applied scientists in space dimensions $d=1,2,3$, in the case of nonlinear parabolic  SPDEs with additive space-time white noise the mathematical literature focuses almost exclusively on the case $d=1$. In fact, in one space dimension, the theory of nonlinear parabolic SPDEs of generic type (\ref{e0}) is well-established in the literature, see e.g.\ \cite{da1992stochastic, walsh1986introduction}. In addition, convergence properties of standard numerical approximation techniques for such equations in $\R^1$ have been thoroughly studied, see \cite{gyongy1995implicit, gyongy1999lattice, lythe2001stochastic, shardlow1999numerical} for finite difference and \cite{walsh2005finite} for finite element methods, respectively. On the other hand, little analytic work has been done on the higher-dimensional cases of these white noise-driven models. With the exception of a few specific cases such as the stochastic Cahn-Hilliard equation (whose well-posedness is established in  \cite{da1996stochastic}), most mathematical studies claim upfront that additive white noise leads to ill-posed equations in $d\geq 2$. Typically, authors then follow one of the following regularization procedures:  they either restrict the analysis to the case of colored noise with a finite spatial correlation length \cite{kohn2007action,erbar2010low,kloeden2010exponential,jentzen2009pathwise,kovacs2010strong}, or they render the equation well-posed by means of the so-called stochastic quantization procedure involving Wick products \cite{da2003strong,chan2000scaling}. 
\\

\begin{figure}[htb!]
   \centering
   \includegraphics[width=13cm]{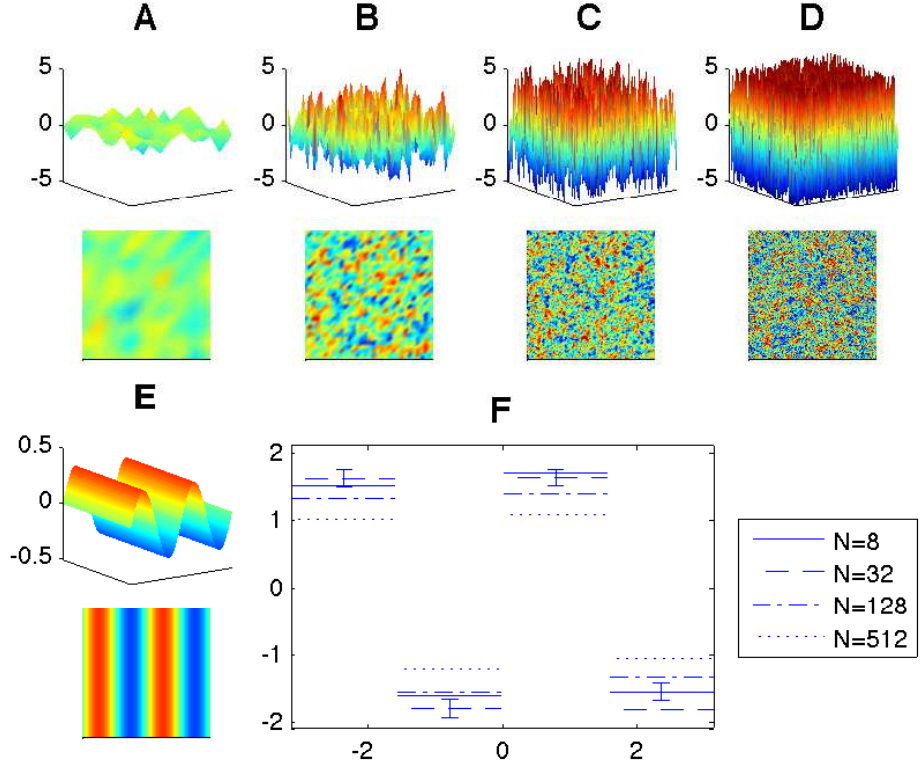} 
   \caption{ {\bf 2D stochastic Allen-Cahn equation}.  {\bf A-D:} Equation (\ref{e45}) is solved over time interval $[0,1]$ using scheme (\ref{e46}) for increasing number of grid points: $N=8$ ({\bf A}), $N=32$ ({\bf B}), $N=128$ ({\bf C}), $N=512$ ({\bf D}). The final fields $u_N(x,y)$ are shown from side and top angles; values outside of $[-5,+5]$ are set to $-5$ and $+5$ respectively. Initial condition $u_0(x,y)=\sin(2x)$; parameters: $\alpha=6.4\cdot 10^{-3}$, $g=0.5$, $\sigma=2^{-3}\,\pi$, $M=10^3$.  {\bf E:} The deterministic Allen-Cahn equation is solved using scheme (\ref{e46}). Initial condition and parameters as above, except  $\sigma=0$, $N=256$. {\bf F:} $[-\pi,+\pi]$ is divided into four subintervals $I_k=[\frac\pi2\,\,(k-1), \frac\pi2 k]$ for $k=-1,\ldots,2$. Simulations A-D are repeated for $120$ realizations of the noise and means of the piecewise constant functions $\tilde{u}_N(z\in I_k):= \int_{I_k} \int_0^{2\pi} u_N(x,y)\,dy\, dx$ are plotted. For each $I_k$, the largest errorbar is shown.} 
   \label{fig:FIG_4}
\end{figure}

\noindent Assembling the above observations, we end up with the following discrepancy: whereas higher-dimensional versions of these nonlinear parabolic SPDE models with additive white noise are commonly used by applied scientists, they are assumed to be ill-posed by the mathematical community.  The goal of this study is to gain a better understanding of what goes wrong in a class of {\it Model A} equations, and we do so by focusing on the specific example of the stochastic Allen-Cahn equation on a bounded domain in $\R^2$, 
\begin{align}\label{e1}
\partial_t u= \Delta u +u - u^3 + \xi,
\end{align}
where $\xi$ is space-time white noise. As discussed in Section \ref{simpli1ana} below, it is unclear how to define a weak solution to this equation. Therefore we take a different approach: we consider a sequence of regularized versions of (\ref{e1}) with corresponding solutions $u_N$, and study the limit $\lim_{N\to\infty} u_N$.  In particular, the regularization procedure is such that the sequence $\set{u_N(t)}_{N=1}^{\infty}$ is $L^2$-valued for all $t\geq0$. As illustrated in Fig.\ \ref{fig:FIG_4}A-D, the sequence of regularizations (shown at time $t=1$) does not seem to converge to a meaningful limit as $N\to\infty$: the underlying deterministic evolution of the field (Fig.\ \ref{fig:FIG_4}E)  gets washed out as $N$ increases and the field becomes highly oscillatory (Fig.\ \ref{fig:FIG_4}C-D). 
Furthermore, the integral of $u_N(t)$ over an arbitrary subset of the spatial domain decays with increasing $N$ (Fig.\ \ref{fig:FIG_4}F). This motivates the following conjecture. 
\begin{conjecture}\label{con1}
The regularizations $u_N$ of the two-dimensional stochastic Allen-Cahn equation (\ref{e1}) converge in probability to the {\bf zero-distribution}, i.e. $\forall t>0$ and all smooth test functions $\phi$:
\begin{align}\label{e1a0}
\lim_{N\to \infty} \ip{u_N(t), \phi} \to 0 \quad \textrm{ in probability,} 
\end{align}
where $\ip{\cdot, \cdot}$ denotes the $L^2$-duality pairing.
\end{conjecture}
\noindent A simple model of this form of convergence consists of the sequence of functions $\set{\sin(Nx)}_{N=1}^\infty$ with  $x\in[0,2\pi]$. All elements of this sequence are smooth functions, but they become highly oscillatory for large $N$, similarly to Fig.\ \ \ref{fig:FIG_4}D. It is easy to show that for all periodic and smooth test functions $\phi$,
\begin{align}\label{e1a01}
\lim_{N\to \infty}\,\, \ip{\sin(N \cdot\,), \phi} = 0,
\end{align}
which means that  $\sin(Nx)$ weakly converges to zero.\\

\noindent The main goal of this study is to provide sound numerical evidence for Conjecture \ref{con1} and to give a heuristic explanation as to why the limit $\lim_{N\to\infty}\,u_N$ is not meaningful.  After a brief review of the deterministic Allen-Cahn equation and the properties of its stochastic version for $d=1$ in Section \ref{example}, we proceed to two simplified versions of the equation in general dimension $d\geq 1$. In Section \ref{simpli1} we revisit the well-studied stochastic heat equation and in Section \ref{simpli2} we consider the decoupled stochastic Allen-Cahn equation - a field of independent, noise-driven particles in double-well potentials. Studying these simplifications provides us with a better understanding of the respective roles played by the nonlinearity and the diffusion in the full stochastic Allen-Cahn equation in two space dimensions. The latter is the subject of Section \ref{AC}, where we provide heuristic arguments and numerical evidence for Conjecture \ref{con1}. In Section \ref{numerics} we present the numerical solution strategies employed in Sections \ref{simpli1}-\ref{AC}. Finally, implications of our work for published numerical studies as well as a presentation of future directions are the subject of Section \ref{conclusions}.


\section{Preliminaries}\label{example}

The Allen-Cahn equation on a domain $\Omega \subset \R^d$,
\begin{align}\label{e1a}
\partial_t u= \Delta u +\frac{1}{\epsilon^2}\ob{u - u^3},
\end{align}
was first introduced by Allen and Cahn to describe a non-conserved order field during anti-phase domain coarsening \cite{allen1979microscopic}. The equation describes the gradient flow of the energy functional
\begin{align}\label{e2}
E[u ]=  \int_\Omega \frac12 \abs{\nabla u}^2 -  \frac{1}{\epsilon^2} V(u) dx
\end{align}
with the particular choice of a double-well potential $V(u)=\frac14 u^4-\frac12 u^2$. The Allen-Cahn equation has been extensively studied in the literature and exhibits a variety of interesting phenomena such as interface motion by mean curvature flow in the limit as $\epsilon \to 0$, see  \cite{evans1992phase}. We are not concerned with such limits and shall hereafter set $\epsilon=1$.\\

\noindent To take into account thermal effects in the gradient flow of (\ref{e2}), a natural approach is to complement the flow with a random forcing term in form of additive space-time white noise $\xi$,
\begin{align}\label{e3}
\partial_t u= \Delta u +u - u^3+ \sigma \xi,
\end{align}
where $\sigma > 0$ is a constant and $\xi$ is a space-time Gaussian random process with mean zero and formal correlation function $\E\cb{\xi(x,t)\xi (x',t')}= \delta(x-x')\,\delta(t-t')$ \cite{shardlow2000stochastic}. Note that in the framework of Hohenberg and Halperin's classification \cite{hohenberg1977theory}, the stochastic Allen-Cahn equation (\ref{e3}) belongs to {\it Model A: systems without conservation laws}.
In the current study, we consider (\ref{e3}) on the torus $\T^d$, i.e. the $d$-dimensional hypercube $[-\pi, +\pi]^d$ subject to periodic boundary conditions. This yields the following initial value problem on $\mathbb{T}^d$:
\begin{align}\label{e4}
\left\{
\begin{array}{rcll}
\partial_t u & =& \Delta u + (u-u^3) + \sigma \xi, & t>0 \\
u &= &u_0,  & t=0.
\end{array}\right.
\end{align}
For the sake of simplicity we assume that $u_0$ is a deterministic function. Later on, the following concise definition of space-time white noise for the torus-setting will be useful. 
\begin{definition}\label{STWN}
Let $T>0$ and consider the measure space $\ob{ \cb{0,T}\times \T^d, \mathcal{B}, \lambda}$, where $\mathcal{B}$ is the Borel $\sigma$-algebra of $\cb{0,T}\times \T^d$ and $\lambda$ the Lebesgue measure. Let $\ob{\Omega, \mathcal{F}, \P}$ be a probability space. A {\bf space-time white noise} on $\cb{0,T}\times \T^d$ is a mapping $\xi: \mathcal B \to L^2(\Omega)$ such that
\begin{enumerate}[(i)]
\item for all $B\in \mathcal B$,  $\xi(B)$ is centred Gaussian with \[ \E\, \cb{\xi(B)}^2 = \lambda(B).\]
\item if $B_1\cap\ldots \cap B_n=\emptyset$, then the $\set{\xi(B_i)}_{i=1}^n$ are independent and \[ \xi \ob{\cup_{i=1}^n B_i} = \sum_{i=1}^n \xi (B_i).\]
\end{enumerate}
\end{definition}
\noindent Recall now that white noise is, formally speaking, the time derivative of the infinite dimensional cylindrical Wiener process $W$ \cite[Sect.\,7.1.2]{peszat2007stochastic}. We can thus rewrite problem (\ref{e4}) as
\begin{align}\label{e4a}
\left\{
\begin{array}{rcll}
du &=& \cb{\Delta u +u -u^3 }dt + \sigma dW, & t>0 \\
u &= &u_0,  & t=0.
\end{array}\right.
\end{align}
This notation will be convenient for our analysis because $W$ admits the spectral decomposition  \cite{da1996stochastic}
\begin{align}\label{e5}
W(t)= \sum_{k\in \Z^n} \beta_k(t) e_k,
\end{align}
where $\set{\beta_k}_{k\in \Z^d}$ are i.i.d.\ Brownian motions and $\set{e_k}_{k\in \Z^d}$ is an orthonormal basis of $L^2(\mathbb{T}^d)$ with respect to the inner product $\ip{h, g}=\int_{[-\pi,+\pi]^d} h \bar{g} dx$, 
\begin{align}\label{e5a}
e_k(x)=\ob{2\pi}^{-d/2}e^{ikx}. 
\end{align}
Before attempting to solve (\ref{e4}), we have to define the notion of a solution. In fact, white noise is too rough to make sense of the equation pointwise, and we proceed formally by integrating (\ref{e4a}) against a smooth test function to obtain the following weak formulation \cite{da1992stochastic}. 
\begin{definition}\label{def1}
Let $\mathcal{H}$ be a Hilbert space. An $\mathcal{H}$-valued process $u(t)$, $t\in \cb{0,T}$, is called a {\bf weak solution} to problem (\ref{e4}) if (i) $\int_0^T \norm{u(t)}dt<+\infty$ for almost all trajectories, and (ii) it satisfies the {\bf weak formulation} 
\begin{align}\label{e7}
\ip{  u(t), \phi } = \ip{ u_0, \phi} + \int_0^t  \cb{\ip{u(s), \Delta \phi}+ \ip{ u(s)-u^3(s), \phi}}ds + \ip{ W(t), \phi},
\end{align}
almost surely, for all $t\in\cb{0,T}$ and for all $\phi\in C^\infty(\T^d)$.  Here, $\ip{\cdot,\cdot}$ denotes the inner product on $\mathcal{H}$.
\end{definition} 
\noindent  For the case $d=1$, existence and uniqueness of the weak solution (\ref{e7}) have been established for $\mathcal{H}=L^2(\T)$ \cite{faris1982large,brassesco1995brownian} (in fact, the solution is almost surely continuous). Therefore, the initial value problem (\ref{e4}) provides a mathematically sound model in one space dimension. However, in higher dimensions ($d\geq 2$) the situation is quite different: the weak solution to the linearized version is only a distribution-valued process, and hence it is unclear whether there is a Hilbert space $\mathcal{H}$ such that a unique weak solution to the nonlinear equation exists. To gain a better understanding of these issues, it is instructive to have a closer look at the stochastic heat equation, a linear version of the full Allen-Cahn model.

\section{Simplified version I: stochastic heat equation}\label{simpli1}

\subsection{Analytic Considerations}\label{simpli1ana}
In this section we focus on the well-studied stochastic heat equation on $\mathbb{T}^d$ \cite{walsh1986introduction},  subject to homogenous initial conditions,
\begin{align}\label{e8}
\left\{
\begin{array}{rcll}
du&=& \cb{\Delta u - u } dt + \sigma dW, & t>0\\
u& = & 0, & t=0.
\end{array}\right.
\end{align}
Using the spectral decomposition of the noise (\ref{e5}), the projection of (\ref{e8}) onto the Fourier modes (\ref{e5a}) yields the following set of stochastic differential equations for $k\in \Z^d$
\begin{align}\label{e9}
\left\{
\begin{array}{rcll}
d\hat{u}_k &=& - \mu_k \hat{u}_k dt + \sigma d\beta_k, & t>0\\
\hat{u}_k& =& 0, & t=0.
\end{array}\right.
\end{align}
where $\mu_k= 1 + \abs{k}^2$ and  $\abs{k}^2=\sum_{j=1}^{d} k_j^2$. In other words, the solution to (\ref{e8}) is represented by a collection of i.i.d.\ Ornstein-Uhlenbeck processes (\cite{gardiner2009stochastic}, p106) whose solutions are given by
\begin{align}\label{e10}
\hat{u}_k(t)=\sigma \int_0^t e^{-\mu_k(t-s)}d\beta_k(s).
\end{align}
These are mean zero Gaussian processes with covariance  ($s>0$)
\begin{align}\label{e11}
\E \, \hat{u}_k(t)\hat{u}_k(t+s) = \frac{\sigma^2}{2\mu_k} e^{-\mu_k s}\cb{1-e^{-2\mu_k t}}.
\end{align}
With this we are now able to calculate the expected value of the $L^2(\T^d)$-norm (denoted by $\norm{\cdot}_0$) of the solution to (\ref{e8}),
\begin{align} \label{e12}
\E \norm{u(t)}_0^2 = \E \sum_{k \in \Z^d} \abs{\hat{u}_k(t)}^2=\sum_{k \in \Z^d} \frac{\sigma ^2 }{2\mu_k} \cb{1-e^{-2\mu_k t}}.
\end{align}
The convergence of the series depends on the summability of $\mu_k^{-1}=(1+\abs{k}^2)^{-1}$ and it is easy to see that $\E \norm{u(t)}^2 < + \infty$ for $d=1$, whereas $\E \norm{u(t)}^2= +\infty$ for $d\geq 2$. In other words, the solution of the one dimensional stochastic heat equation (\ref{e8}) takes values in $L^2(\T^1)$ almost surely -- a result that does not hold true in higher dimensions. \\

\noindent Let us now have a closer look at the case $d\geq 2$, and more precisely at the rate of divergence of the sum in (\ref{e12}). For fixed $t>0$,
\begin{align}\label{e13}
\E\norm{u(t)}_0^2 & \sim \sum_{\abs{k}\leq N} \frac{1}{1+\abs{k}^2}  \sim\int_0^N \frac{r^{d-1}}{1+r^2}  dr,
\end{align}
where we use the following definition: $f_N\sim g_N$ if there exist two constants $c, C>0$ such that $c \leq f_N/g_N \leq C$
 for $N$ sufficiently large. (\ref{e13}) implies that the divergence is logarithmic for $d=2$ and polynomial for $d\geq 3$. But even though it is not $L^2(\T^d)$-valued, $u(t)$ might be well-defined in a larger space. In this sense, a natural enlargement of $L^2(\T^d)$ is given by the Sobolev spaces of negative order, $H^s(\T^d)$ for $s<0$ (e.g.\ \cite{dautray2000mathematical}, p96).
\begin{definition}[Sobolev Spaces]
Let $s\in \mathbb{R}$. Then the Sobolev space $H^s(\T^d)$ is defined by means of a weighted norm as
\begin{align}\label{e14} 
H^s(\mathbb{T}^d):= \set{ f: \mathbb{T}^d \to \R \quad : \quad \norm{f}_s^2=\sum_{k\in \Z^d} \ob{1+\abs{k}^2}^s \abs{\hat{f}_k}^2 <+\infty},
\end{align}
where  $\hat{f}_k$ are the generalized Fourier coefficients. In particular, $L^2(\T^d)=H^0(\T^d)$. 
\end{definition} 
\noindent Revisiting the calculation (\ref{e13}) in $H^s(\T^d)$ yields 
\begin{align}\label{e15}
\E \norm{u(t)}_s^2 \sim \int_0^N \frac{r^{d-1}}{(1+r^2)^{1-s}}dr<+\infty, \qquad \forall s<1-d/2,
\end{align}
i.e. the solution to the $d$-dimensional version of (\ref{e8}) takes values in $H^s(\T^d)$ almost surely, for all $s<1-d/2$.
We summarize these estimates in the following result \cite{walsh1986introduction}.
\begin{theorem}\label{prop1}
Let $d\geq1$. For all $s<1-d/2$, the solution to the $d$-dimensional stochastic heat equation takes values in $H^s(\T^d)$ almost surely. In particular, $d=2$ is the borderline case:  the $L^2$-divergence in (\ref{e13}) is logarithmic, and for all $t>0$ we find that $u(t)\in H^s(\T^d)$ almost surely, for all $s<0$. 
\end{theorem}
\noindent Theorem \ref{prop1} illustrates the smoothing properties of the heat kernel: whereas space-time white noise only takes values in $H^s(\T^d)$ for $s<-d/2$, the action of the heat semigroup improves the regularity such that $u(t)\in H^s(\T^d)$ for  $s<-d/2+1$. In view of the nonlinear stochastic Allen-Cahn equation (\ref{e4}) in dimensions $d\geq2$, the relevance of Proposition \ref{prop1} is the following: for negative $s$,  $H^s(\T^d)$ is a space of distributions, and it is in general not possible to multiply two distributions in a meaningful way \cite{schwartz1954impossibilite}. In consequence, taking the cube of the linear solution does not make sense and we anticipate the nonlinear equation  to be ill-posed. This observation provides a heuristic argument as to why white noise-driven nonlinear parabolic SPDEs in higher dimensions are generally suspected to be ill-posed. However, it is not a rigorous proof and the question remains as to whether it is possible to find a space in which there exists a weak solution to the stochastic Allen-Cahn equation in higher dimensions.

\subsection{Simulations}\label{simpli1num}
\begin{figure}[htb!]
   \centering
   \includegraphics[width=10cm]{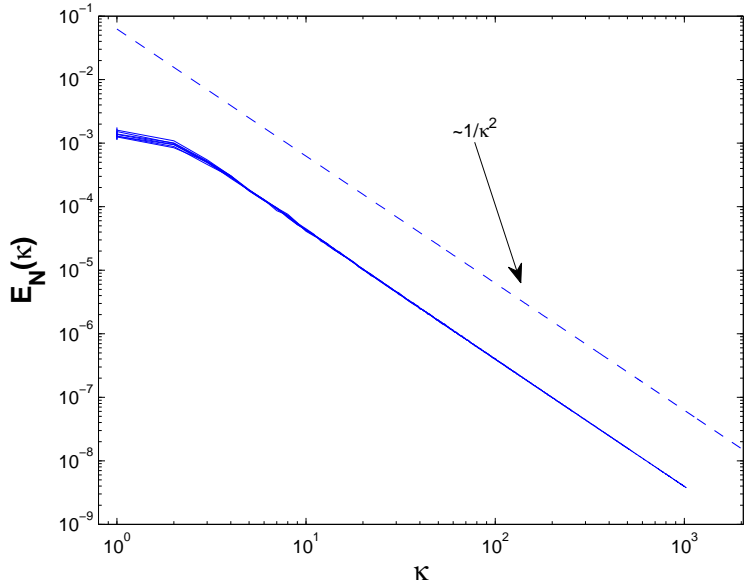} 
   \caption{ {\bf 2D stochastic heat equation}. Equation (\ref{e40}) is solved using scheme (\ref{e42}). Parameters: $g=1$, $\alpha=0.5$, $\sigma=\pi/50$, $T=1$, $M=2\cdot 10^3$. For each $N=2^n$ ($n=5,\ldots,11$), the sample mean of $E_N(\kappa)$ over $40$ simulations is plotted for $\kappa=1,\ldots,N/2$. Largest error bar is shown for each $N$, and we observe convergence within errorbars. Dotted line is of slope $1/\kappa^2$.} 
   \label{fig:FIG_1}
\end{figure}

To illustrate the results of the previous section, and to facilitate comparisons with subsequent results, we present now a numerical experiment on the two-dimensional stochastic heat equation (\ref{e8}). We discretize the periodic domain $\T^2$ using $N^2$ grid points $\{x_j=\frac{2\pi}{N} \,j\, : \ j=(j_1,j_2),\,\, j_i=-N/2,\ldots, N/2-1\}$ and denote by $u_N(x_j,t)$ the numerical approximation at grid-point $x_j$ and time $t$. Since we are mostly interested in the spatial regularity of the solution, we fix $T>0$, integrate the equation numerically over $[0,T]$, and denote $u_N(x_j)\equiv u_N(x_j,T)$. To estimate the regularity of the approximation $u_N$, we first need to define the finite-dimensional analogue of the Sobolev norm (\ref{e14}). To do this, we introduce the discrete Fourier transform of $u_N$ as 
\begin{align}\label{e15a}
\hat{u}_N(k)= \frac{1}{N} \sum_{j_1=-N/2}^{N/2-1}\,\sum_{j_2=-N/2}^{N/2-1} u_N(x_j) \exp{\ob{\frac{2\pi i}{N} k\cdot j}},
\end{align} 
and the discrete inverse Fourier transform as
\begin{align}\label{e15b} 
u_N(x_j)= \frac{1}{N} \sum_{k_1=-N/2}^{N/2-1} \,\sum_{k_2=-N/2}^{N/2-1}\hat{u}_N(k) \exp{\ob{-\frac{2\pi i}{N} k\cdot j}},
\end{align}
where $\set{ k=(k_1,k_2) \, :\, k_i=-N/2,\ldots,N/2-1}$ are the wave vectors. We can now define the  $H^s(\T^d)$-norm for discrete functions defined on the $N^2$ grid points. 
\begin{definition}[Discrete Sobolev Norm]\label{defnorm}
Let $N\geq2$ be an even integer. Let $u_N$ be a discrete function defined on the $N^2$ grid points of $\T^2$. The discrete $H^s(\T^2)$-norm of $u_N$ is defined as 
\begin{align}\label{e15c}
\dnorm{u_N}_s^2:= \rho^{-1} \sum_{k_1=-N/2}^{N/2-1} \,\sum_{k_2=-N/2}^{N/2-1} \ob{1+|k|^2}^s\abs{\hat{u}_N(k)}^2,
\end{align}
where $\rho^{-1/2}:=2\pi/N$ is the grid spacing.
\end{definition}
\begin{remark} To see that $\rho^{-1}$ is the correct scaling in (\ref{e15c}), consider the case $s=0$. Using Parseval's equality we find 
\begin{align}\label{e15a1}
\dnorm{u_N}_0^2=& \,\rho^{-1} \sum_{k_1=-N/2}^{N/2-1}\,\sum_{k_2=-N/2}^{N/2-1}\abs{\hat{u}_N(k)}^2=  \rho^{-1}  \sum_{x_1=-N/2}^{N/2-1} \,\sum_{x_2=-N/2}^{N/2-1}\abs{u_N(x)}^2\nonumber \\ 
&\xrightarrow{N\to\infty}  \int_{[0,2\pi]^2} \abs{u}^2 dx = \norm{u}_0^2.
\end{align} 
\end{remark}
\noindent The following notion of the {\it radial energy density} will be useful for the graphical illustration of the numerics.
\begin{definition}[Radial Energy Density in Fourier Space]
Let $N\geq2$ be an even integer and let $u_N$ be a discrete function defined on the $N^2$ grid points of  $\T^2$. The radial energy density in Fourier space, $E_N$, is defined as
\begin{align}\label{e15a2}
E_N(\kappa): = \frac{\rho ^{-1}}{\abs{A_\kappa}} \sum_{k\in A_\kappa} \abs{\hat{u}_N(k)}^2, \qquad \kappa=1,\ldots,N/2-1,
\end{align}
where $A_\kappa=\set{ x \in \R^2 : (\kappa-1)^2\leq \abs{x}^2 \leq \kappa^2}$ is the $\kappa$-th annulus in $\R^2$ and $\abs{A_\kappa}$ is the cardinality of $\Z^2\cap A_\kappa$. 
\end{definition}
\noindent Using the energy density (\ref{e15a2}), we can replace the $H^s(\T^2)$-norm (\ref{e15c}) with an equivalent discrete norm defined as
\begin{align}\label{e15a3}
\dnorm{u_N}_s^2 =\sum_{\kappa=1}^{N/2-1} \kappa\, E_N(\kappa) (1+\kappa^2)^s.
\end{align}
From this we conclude that in the continuum limit as $N\to\infty$, the regularity of $u_N$ is determined by the decay of the radial energy density $E_N$. In Figure \ref{fig:FIG_1}, the density $E_N$ for the corresponding numerical solution $u_N$ of the stochastic heat equation is plotted for increasing values of $N$: we observe convergence within the errorbars, and in particular the $E_N(\kappa)$ decay like $1/\kappa^2$ for $\kappa$ sufficiently large.
By (\ref{e15a3}), this suggests that $\lim_{N\to \infty} \E \norm{u_N}_0^2 = +\infty$, whereas $\lim_{N\to\infty} u_N \in H^s(\T^2)$ almost surely, $\forall s<0$. Note that these limits are in perfect agreement with Theorem \ref{prop1}. A detailed discussion of the employed numerical scheme is presented in Section \ref{numsimpli1}.

\section{Simplified version II: decoupled stochastic Allen-Cahn equation}\label{simpli2}
\subsection{Analytic Considerations}\label{simpli2ana}

In this section we study a second simplification of the stochastic Allen-Cahn equation in $\R^d$. This time, we keep the nonlinearity but neglect the diffusion term, i.e. we consider 
\begin{align}\label{e16}
du= \cb{u-u^3}dt + \sigma dW,
\end{align}
on the domain $\T^d$ and subject to homogenous initial conditions.  Even though this equation does not provide a meaningful physical model, it is instructive in view of the discussion of the full stochastic Allen-Cahn equation in Section \ref{AC}. In fact, we will see that approximations to the stochastic Allen-Cahn equation and its decoupled version (\ref{e16}) share similar characteristics. However, the latter is more tractable due to the spatial decoupling: at each point in space there is a particle, driven by a Brownian motion and confined by a double-well potential. Note that the particles are i.i.d.\ because the collection of driving Brownian motions arises from space-time white noise. It is not clear if there exists a solution space which allows for the definition of a weak solution to (\ref{e16}). We circumvent this difficulty by first discretizing the equation in direct space and then passing to the continuum limit. We choose an even integer $N$, introduce an equidistant mesh of size $\rho^{-1/2}:=\Delta x=2\pi/N$ along each dimension of $\T^d$, and consider the finite-dimensional approximation $u_N$, which is defined on the $N^d$ grid points $\set{x_j}$. At each grid point, the evolution equation is
\begin{align}\label{e17}
du_N(x_j)= -V'\ob{u_N(x_j)}dt + \sigma \rho^{d/4} d\beta_j,
\end{align}
where $V(x):=\frac14 x^4-\frac12 x^2$ is a double-well potential, $\set{\beta_j}$ is a collection of i.i.d. standard Brownian motions. The scaling $\rho^{d/4}$ is due to the spatial discretization of white noise. The mesh partitions the domain into $N^d$ hyper-cubes of volume $\rho^{-d/2}$, and from Definition \ref{STWN} it follows that the average noise on each hyper-cube (at fixed time $t$) is distributed as $ \rho^{d/2}\mathcal{N} (0,\rho^{-d/2}) \sim \rho^{d/4} \mathcal{N} (0,1)$.  As the number of grid points increases, the deterministic term in (\ref{e17})  remains unaltered while the noise intensity increases.  Therefore, the variance of the solution $u_N(t)$ is unbounded in the limit as  $N\to\infty$. Unless the potential $V$ is quadratic, it is not possible to write down the dynamic solution to (\ref{e17}) in closed form; instead we focus on the stationary solution (as $t\to\infty$) whose probability distribution function is given by 
\begin{align}\label{e18}
p_N(y)=\frac{1}{\mathcal{N}} \exp{\set{-\frac{2V(y)}{\sigma^2\rho^{d/2}}}},
\end{align}
with $\mathcal{N}$ the normalization constant \cite{gardiner2009stochastic}. We prove now the following result.
\begin{theorem}\label{prop3}
Let $d\geq 1$ and $V(x)=\frac14x^4-\frac12 x^2$. Then there exists a constant $K>0$ such that  the stationary solution $u_N$ of the regularized problem (\ref{e17})  admits the following limits
\begin{align}\label{e24}
\lim_{N\to \infty} \E\, \dnorm{u_N}_s^2 = \left\{ 
\begin{array}{cl} 
0 & \text{if} \quad s< -d/4\\
K<+\infty & \text{if} \quad s=-d/4\\
+\infty &\text{if} \quad s> -d/4,
\end{array} \right.
\end{align}
where $\dnorm{\cdot}_s$ is the discrete Sobolev norm of Definition \ref{defnorm}, generalized to $d$ dimensions (\ref{e23}).
\end{theorem}
\begin{proof} The main part of the proof is established in Lemma \ref{appA} in the Appendix. We first generalize the discrete Fourier transform (\ref{e15a}) and its inverse (\ref{e15b}) to $d$ dimensions:
\begin{align}\label{e24a}
\hat{u}_N(k)= \frac{1}{N^{d/2}} \sum_{j_i=-N/2}^{N/2-1} u_N(x_j) \exp{\ob{\frac{2\pi i}{N} k\cdot j}},
\end{align} 
and the inverse transform
\begin{align}\label{e24b} 
u_N(x_j)= \frac{1}{N^{d/2}} \sum_{k_i=-N/2}^{N/2-1} \hat{u}_N(k) \exp{\ob{-\frac{2\pi i}{N} k\cdot j}}.
\end{align}
Using (\ref{e24a}) together with point (ii) of Lemma \ref{appA}, we find that the stationary solution $u_N$ of (\ref{e17}) satisfies the following properties:
\begin{align}
 \mathbb{E}\,\hat{u}_N(k)&=0,\quad\forall k=-N/2, \ldots, N/2-1, \quad \forall N\geq 2 \label{e23a}\\
\mathbb{E}|\hat{u}_N(k)|^2& \sim  \sigma N^{d/2}, \quad \forall k=-N/2, \ldots, N/2-1, \quad \forall N\geq 2. \label{e23b}
\end{align}
From this we see that as we refine the grid, $\lim_{N\to\infty} \E |\hat{u}_N(k)|^2=\infty$ for all $k$ and $d\geq 1$. Generalizing the discrete version of the $H^s(\T^d)$-norm (\ref{e15c}) to general dimension $d\geq1$, and using the growth rate (\ref{e23b}), we find
\begin{align}\label{e23}
\E \dnorm{u_N}_s^2 &=  \rho^{-d/2} \sum_{k_i=-N/2}^{N/2-1} \ob{1+|k|^2}^s\abs{\hat{u}_N(k)}^2 \nonumber \\ &\sim \rho^{-d/4} \int_0^N \ob{1+r^2}^s r^{d-1}dr \sim N^{2s+d/2}.
\end{align}
The limits in (\ref{e24}) follow now easily.
\end{proof}
\noindent In particular, (\ref{e24}) implies that the continuum limit of $u_N$ does not take values in $H^s(\T^d)$ for any $s>-d/4$, and in particular is not in $L^2(\T^d)$. 
Furthermore, recalling that the dual of  $H^s(\T^d)$ is $H^{-s}(\T^d)$, and that $C^\infty(\T^d)\subset H^s(\T^d)$ $\forall s \in \R$, it follows immediately that 
\begin{corollary}\label{cor4}
Let $d\geq 1$ and $V(x)=\frac14x^4-\frac12x^2$. Then the stationary solution $u_N$ of the regularized problem (\ref{e17}) converges in probability to the zero-distribution in the sense of the definition given in Conjecture \ref{con1}.
\end{corollary}

\subsection{Simulations}\label{simpli2num}
\begin{figure}[htb!]
   \centering
   \includegraphics[width=10cm]{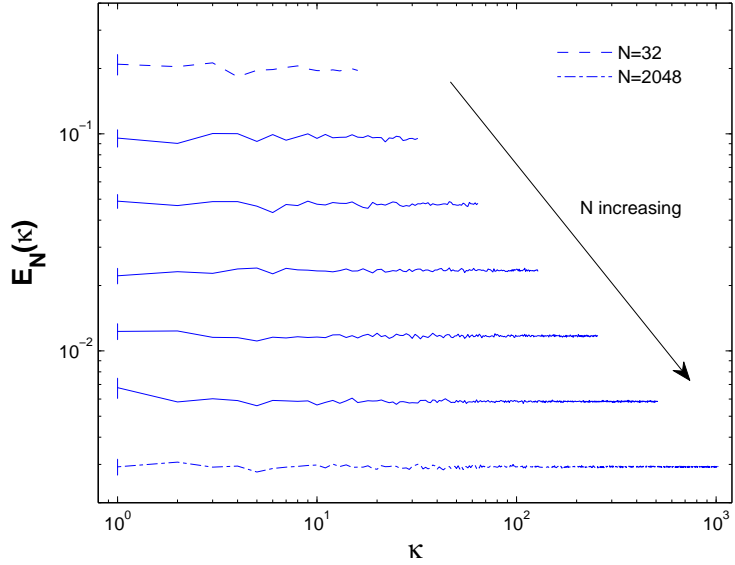} 
   \caption{ {\bf 2D decoupled stochastic Allen-Cahn equation }. Equation (\ref{e43}) is solved using scheme (\ref{e44}). Parameter values: $g=0.1$, $\sigma=\pi/5$, $T=2$, $M=4\cdot 10^3$. For each $N=2^n$ ($n=5,\ldots,11$), the sample mean of $E_N(\kappa)$ over $40$ simulations is plotted for $\kappa=1,\ldots,N/2$. Largest error bar is shown for each $N$. }
   \label{fig:FIG_2}
\end{figure}

Similarly to the case of the stochastic heat equation in Section \ref{simpli1num}, we compare our results on the decoupled stochastic Allen-Cahn equation now to numerical experiments. According to (\ref{e23b}) the mean of the radial energy density $E_N$ (\ref{e15a2}) decays like
\begin{align}\label{e24c}
\E\, E_N(\kappa) \sim  \frac1N, \qquad \kappa=1,\ldots,N/2-1, \qquad N\to \infty.
\end{align}
independently of the spectral radius $\kappa$. This is illustrated by the simulations in Fig.\ \ref{fig:FIG_2}, where the energy density is plotted for various values of $N$.
Taking into account the loglog-scale of the plot, the regular spacing of the $E_N$ for increasing $N$ ($N=2^j$ for $j=5,\ldots,11$) illustrates the decay in (\ref{e24c}). A detailed discussion of the employed numerical scheme is presented in Section \ref{numsimpli2}.

\section{2D Stochastic Allen-Cahn equation}\label{AC}

\subsection{Heuristic Considerations}\label{ACana}
We turn our attention to full the stochastic Allen-Cahn equation in two dimensions, i.e. we consider 
\begin{align}\label{e25}
du= \cb{\Delta u +u -u^3}dt +\sigma dW
\end{align}
on $\T^2$, subject to homogenous initial conditions. As outlined at the end of Section \ref{simpli1}, it is unclear how to make sense of this equation: in fact, we suspect the continuum version of (\ref{e25}) to
be ill-posed with respect to any reasonable definition of a weak solution. Therefore, we take a different approach and tackle the problem from the perspective of numerical analysis. By cutting off high frequency modes in the noise,
\begin{align}\label{e26}
dW_N= \sum_{\abs{k}\leq N} \beta_k(t) e_k,
\end{align}
we obtain the following family of regularized problems
\begin{align}\label{e27}
du_N= \cb{\Delta u_N +u_N -u_N^3}dt +\sigma dW_N.
\end{align}
For fixed $N$, the solution to (\ref{e27}) exists and is unique, but what can we say about the limit of $u_N$ as $N\to \infty\,$? In the remainder of this section, we discuss this limit and provide evidence for Conjecture \ref{con1}, the main subject of this study. The arguments hereafter will invoke elements of the theory of stochastic quantization and we refer the reader to \cite{da2007wick} for a rigorous introduction to the formalism.\\

\noindent At the end of Section \ref{simpli1} we explained that the term $u_N^3$ in (\ref{e27}) was expected to be undefined in the limit $N\to\infty$. To circumvent this issue, we introduce a constant $C_N>1$, depending only on $N$, and rewrite equation (\ref{e27}) as
\begin{align}\label{e28}
du_N= \cb{\Delta u_N - \ob{C_N -1} u_N -u_N\ob{u_N^2-C_N}}dt +\sigma dW_N.
\end{align}
The key idea is to choose the constant $C_N$ wisely, so that the term $u_N\ob{u_N^2-C_N}$ approaches a well defined limit. More precisely, we choose $C_N$ such that $u_N\ob{u_N^2-C_N}$ converges to $:u^3:$, the renormalized cube of $u$, which we briefly explain here. Given the stationary measure $\mu_N$ of the underlying linear equation
\begin{align}\label{e29}
dv_N= \cb{\Delta v_N - \ob{C_N -1} v_N}dt +\sigma dW_N,
\end{align}
$x\mapsto :x^n:$  denotes the $n$-th Wick power with respect to the limit measure $\lim_{N\to\infty} \mu_N$.
In contrast to $\lim_{N\to\infty}u_N^3$, the renormalized cube $\lim_{N \to\infty} u_N\ob{u_N^2-C_N}=:u^3:$ is well-behaved in the sense that its generalized Fourier coefficients $\ip{:u^3:,e_k}$ are uniformly $L^2$-bounded in the probabilistic sense. But what is the correct choice for $C_N$? In fact, its job is to subtract the diverging parts that cause the unmodified cube to be ill-posed, and it is given by
\begin{align}\label{e29a}
C_N=3 \E_{\mu_N} \abs{u_N}^2,
\end{align} 
where $\E_{\mu_N}$ denotes the expected value with respect to $\mu_N$. To determine the measure $\mu_N$, we project (\ref{e29}) onto the Fourier basis to obtain a collection of i.i.d.\ Ornstein-Uhlenbeck processes. From this it easily follows that $\mu_N$ is a mean zero Gaussian measure and according to (\ref{e11}) its covariance operator is given by
\begin{align}\label{e29b}
\hat{C}_N e_k= \left\{ 
\begin{array}{cl} 
\frac{\sigma^2}{2 (C_N -1 +\abs{k}^2)} e_k & \text{if  } \abs{k}\leq N,\\ 0 &\text{if  } \abs{k}>N.
\end{array} \right.
\end{align}
Combining (\ref{e29a}) and (\ref{e29b}) yields the following equation for $C_N$ 
\begin{align}\label{e30}
C_N= \frac{3}{8\pi^2}  \sum_{\abs{k}\leq N} \frac{\sigma^2}{ C_N-1+\abs{k}^2}.
\end{align}
For all $N\geq2$, equation (\ref{e30}) admits a unique positive fixed point, $C_N\to \infty$ as $N \to \infty$, and the rate of divergence is logarithmic 
\begin{align}\label{e31}
C_N =\frac{3\sigma^2}{4\pi}  \log{N}+o(1).
\end{align}
To show (\ref{e31}), rewrite the sum in (\ref{e30}) as a Riemann sum, estimate the latter by an integral approximation (for large $N$), and evaluate the integral to get
\begin{align}\label{e31a}
C_N \sim \frac{3\sigma^2}{8\pi} \log{\ob{\frac{N^2}{C_N}}}=\frac{3\sigma^2}{4\pi}\cb{\log{N}-\frac12 \log{C_N}}.
\end{align}
Making now the Ansatz $C_N=\Lambda_N \log{N}$ and plugging it into (\ref{e31a}) yields the estimate (\ref{e31}). With this choice of $C_N$, we go back to the original equation (\ref{e28}) and project it onto the Fourier basis to obtain
\begin{align}\label{e32}
d\hat{u}_N(k)= -\ob{C_N-1+\abs{k}^2}\hat{u}_N(k) dt - \ip{ u_N (u_N^2-C_N), e_k} dt + \sigma_N(k) d\beta_k(t),
\end{align}
where $\sigma_N(k):=\sigma$ if $\abs{k}\leq N$ and $\sigma_N(k):=0$ if $\abs{k}>N$. 
Recalling that $C_N \to \infty$ and $  \ip{ u_N (u_N^2-C_N), e_k}\to \ip{ :u^3:, e_k}$ (which are uniformly $L^2(P)$-bounded), the first term in (\ref{e32}) dominates the second term as $N$ becomes large. Hence
\begin{align}\label{e33}
d\hat{u}_N(k)\sim -\ob{C_N-1+ \abs{k}^2} \hat{u}_N(k) dt + \sigma_N(k) d\beta_k(t).
\end{align}
The solution to (\ref{e33}) is an Ornstein-Uhlenbeck process for $\abs{k}\leq N$ and a decaying exponential for $\abs{k}>N$, and we conclude that for large times 
\begin{align}\label{e34}
\E \abs{\hat{u}_N(k)}^2 \sim \frac{\sigma^2_N(k)}{2\ob{C_N-1+\abs{k}^2}}, \qquad t\to \infty.
\end{align}
From this we deduce two important facts: first, for a fixed cut-off $N$, the energy in Fourier space essentially decays as
\begin{align}\label{e35}
\E \abs{\hat{u}_N(k)}^2 \sim \frac{1}{\abs{k}^2}, \qquad N\gg 1 \text{ fixed, } \abs{k}\leq N.
\end{align}
Second, for a given mode $k$, the decay as $N \to\infty$ goes as 
\begin{align}\label{e36}
\E \abs{\hat{u}_N(k)}^2 \sim \frac{1}{\log{N}}, \qquad  \text{ fixed } k.
\end{align}
Furthermore, given the joint dependence on $k$ and $N$ in (\ref{e34}), we can estimate the continuum limit of the $H^s(\T^2)$-norm of $u_N$. For $s=0$ we get the $L^2(\T^2)$-norm
\begin{align}\label{e37}
\E\norm{u_N}_0^2=\frac{\sigma^2}{2} \sum_{\abs{k}\leq N} \E\abs{\hat{u}_N(k)}^2 \sim \frac{\sigma^2\pi}{2} \log{\ob{1+\frac{N}{\sqrt{C_N}}}} \to \infty \quad \text{as } N\to \infty,
\end{align}
because (\ref{e31}) implies that $N^{-1} \sqrt{C_N}\to 0$ as $N\to\infty$. On the other hand, if $s<0$, we can use Lebesgue's dominated convergence theorem to get
\begin{align}\label{e38}
\E\norm{u_N}_s^2 =& \,\,\frac{\sigma^2}{2} \sum_{\abs{k}\leq N} \ob{1+\abs{k}^2}^s \E\abs{\hat{u}_N(k)}^2 \nonumber \\
\sim &\,\,\sigma^2 \pi \int_0^\infty \mathbf{1}_{[0,N]}\,\, \frac{x^{1+2s}}{x^2+C_N} dx \to 0 \quad \text{as } N\to\infty, \quad \forall s<0.
\end{align}
Note that the results (\ref{e37}) and (\ref{e38}) have been derived under the assumption that the field has reached the steady-state distribution ($t\to \infty$). However, since the decay rate $(C_N -1 + \abs{k}^2)$ in (\ref{e33}) tends to infinity in the continuum limit, the damping occurs infinitely fast as $N \to \infty$, and we anticipate the limits (\ref{e37}) and (\ref{e38}) to hold for any finite time $t>0$. In other words we propose the following conjecture.
\begin{conjecture}\label{con2}
The solution $u_N$ to the regularized problem (\ref{e27}) admits the following limit. For all $t>0$, 
\begin{align}\label{e39}
\lim_{N \to \infty} \E \norm{u_N(t)}_s^2 = \left\{ 
\begin{array}{cl} 
+\infty & \text{if} \quad s\geq 0, \\
0 & \text{if} \quad s<0.\\
\end{array} \right. 
\end{align}\end{conjecture}
\noindent An immediate corollary to Conjecture \ref{con2} is Conjecture \ref{con1}: $u_N$ converges in probability to the zero-distribution. Even though these two conjectures do not tell us anything about the nature of the continuum solution of the stochastic Allen-Cahn equation (\ref{e25}), they suggest that the approximation sequence (\ref{e27}) does not admit a meaningful limit: even though the limit takes values in the distributional spaces $H^s(\T^2)$ almost surely, for all $s<0$, it is a trivial limit as it projects every test function onto $0$. \\

\noindent To gain a better understanding of the nature of this pathology, it is instructive to have a closer look at the roles played by the three main ingredients of the regularized stochastic Allen-Cahn equation (\ref{e27}): diffusion, nonlinearity, and noise. Considering only the noise, the solution to the corresponding equation, $du_N=dW_N$, converges to the cylindrical Wiener process $W$ as $N\to\infty$. The latter takes values in $H^s(\T^2)$ almost surely, for all $s<-1$.  Adding the nonlinearity to the noise, we obtain the decoupled  equation $du_N=\cb{u_N-u_N^3}dt + dW_N$ of Section \ref{simpli2}. In this case, the $u_N$ converge in probability to a distribution $u_c$ which is more regular than $W$: $u_c\in H^s(\T^d)$ almost surely, for all  $s\leq -1/2$; however, $u_c$ is the zero-distribution. Finally, when adding the Laplacian to obtain the regularized stochastic Allen-Cahn equation,  $du_N=\cb{\Delta u_N +u_N-u_N^3}dt + dW_N$, the regularity of the limit as $N\to\infty$ is conjectured to improve: (\ref{e39}) implies that $u_N$ converges in probability to a distribution $u_c$ with values in $H^s(\T^2)$ almost surely, for all $s<0$. At the same time, the presence of the Laplacian is not sufficient to render the limit meaningful: $u_c$ is again the zero-distribution.

\subsection{Simulations}\label{ACnum}
\begin{figure}[h!]
   \centering
   \includegraphics[height=10cm]{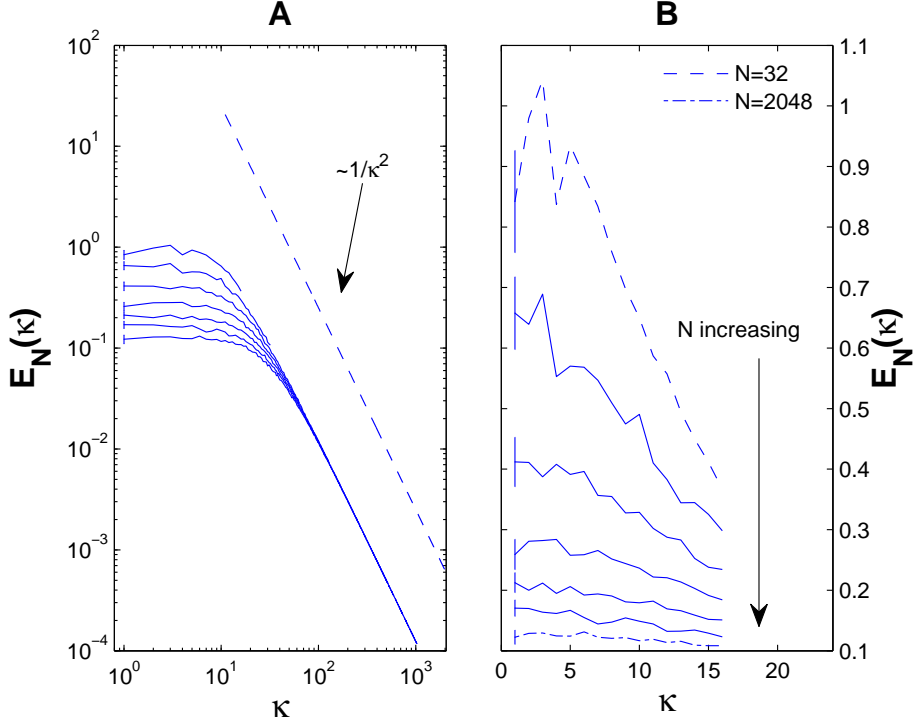} 
   \caption{ {\bf 2D stochastic Allen-Cahn equation}. The stochastic Allen-Cahn equation (\ref{e45}) is solved in two space dimensions using scheme (\ref{e46}). Parameter values: $\alpha=6.4\cdot 10^{-3}$, $g=0.5$, $\eps=2\pi/5$, $T=1$, $M=2\cdot 10^3$. {\bf A. } For each $N=2^n$ ($n=5,\ldots,11$), the sample mean of $E_N(\kappa)$ over $40$ simulations is plotted for $\kappa=1,\ldots,N/2$. Largest error bar is shown for each $N$. {\bf B.} We zoom in on the first 16 modes of panel A. For each $N=2^n$ ($n=5,\ldots,11$), the sample mean of $E_N(\kappa)$ over $40$ simulations is plotted for $\kappa=1,\ldots,16$. Largest error bar is shown for each $N$. } 
   \label{fig:FIG_3}
\end{figure}

The results of the previous section are based on heuristic arguments in absence of rigorous proofs. To put them on more solid ground, we performed a series of numerical experiments. 
 Similarly to the numerical analysis for the simplified equations in Sections \ref{simpli1num} and \ref{simpli2num}, we integrated the stochastic Allen-Cahn equation over the time interval $[0,1]$ and computed the radial energy density $E_N$ at $t=1$, see Figure \ref{fig:FIG_3}. 
A first observation concerns the decay of $E_N(\kappa)$ for fixed $N$: Figure \ref{fig:FIG_3}A shows that this decay is proportional to $\sim 1/\kappa^2$ for large $\kappa$, which is in agreement with the estimate (\ref{e35}). A second observation concerns the decay of $E_N(\kappa)$ for fixed $\kappa$. In Fig.\  \ref{fig:FIG_3}B we zoom into the interval $\kappa\leq 16$ of Fig.\  \ref{fig:FIG_3}A: rather than converging to a finite value, the energy density slowly decays with increasing $N$. In agreement with estimate (\ref{e36}) this decay is logarithmically slow, and hence the expected convergence to zero is not visible within the range of computationally tractable $N$. A detailed discussion of the employed numerical scheme is presented in Section \ref{numAC}. \\

\noindent A comparison between the three figures -- Fig.\ \ref{fig:FIG_1} for the stochastic heat equation, Fig.\ \ref{fig:FIG_2} for the decoupled stochastic Allen-Cahn equation and Fig.\ \ref{fig:FIG_3} for the stochastic Allen-Cahn equation -- suggests that the latter inherits and combines the main features of the two simplified versions: the $\sim 1/\kappa^2$ decay of the Fourier spectrum as dictated by the Laplacian, and the decay (to zero) of individual Fourier modes as dictated by the nonlinearity.

\section{Numerics}\label{numerics}
To perform the two-dimensional numerical experiments presented in Sections \ref{simpli1num}, \ref{simpli2num} and \ref{ACnum}, the periodic domain $[-\pi, +\pi]^2$ was discretized into a regular grid of $N^2$ points ($N$ even), and the interval of integration $[0,T]$ was divided into $M$ subintervals of length $\Delta t=T/M$. Depending on the nature of the equation, the time-stepping was performed in either direct space for the decoupled Allen-Cahn equation (\ref{e16}), in discrete Fourier space for the heat equation (\ref{e8}), or using both spaces for the Allen-Cahn equation (\ref{e25}). Note that for all presented experiments we had verified convergence in time by direct comparison of the solutions obtained with time steps $2\Delta t$ and $\Delta t$, respectively. Before presenting the numerical schemes, a comment about diffusion and reaction rates is in order. Note that equations (\ref{e8}), (\ref{e16}) and (\ref{e25}) are free of constants, except for the noise intensity $\sigma$. In fact, we had eliminated all other constants by rescaling time and space in the corresponding evolution equations. For numerical purposes however, it is convenient to use the same space and time scales throughout the experiments, and hence we re-introduced diffusion and reaction rates (as specified below) in the previously rescaled terms. 

\subsection{Stochastic Heat Equation}\label{numsimpli1}
Introducing a diffusion constant $\alpha>0$ and a decay rate $g>0$, the stochastic heat equation (\ref{e8}) becomes 
\begin{align}\label{e40}
du=\cb{\alpha \Delta u- g u} dt + \sigma dW.
\end{align}
Linearity of the equation allows us to project it onto the $N^2$ modes in discrete Fourier space as 
\begin{align}\label{e41}
d\hat{u}_N(k)= - (g+\alpha \abs{k}^2) \hat{u}_N(k) dt + \sigma d\hat{W}_N(k), 
\end{align}
where $W_N$ is defined in (\ref{e26}) and $k=(k_1,k_2)$. An implicit-explicit scheme was used for the time-stepping, combining the explicit Euler-Maryama scheme \cite{higham2001algorithmic} for the noise with the implicit trapezoidal scheme for the reaction and diffusion terms. With the notations $\hat{u}^m_N(k):=\hat{u}_N(k, t_m)$ and $\Delta x:= 2\pi/ N$,  the scheme reads
\begin{align}\label{e42}
\hat{u}^{m+1}_N(k)=\hat{u}^{m}_N(k) - \frac{\Delta t}{2} (g+\alpha \abs{k}^2)\cb{\hat{u}^{m+1}_N(k)+\hat{u}^{m}_N(k)}+\sigma \frac{\sqrt{\Delta t} }{\Delta x} \xi_k^m,
\end{align}
where $\xi_k^m=\frac{1}{\sqrt{2}}\ob{\eta_k^{m}+ i\, \zeta_k^{m}}$, and $\set{\eta_k^{m}}$ and $\set{\zeta_k^m}$ are real i.i.d.\ $\mathcal{N}(0,1)$ random variables. This representation of white noise on the grid is dictated by the choice of complex eigenfunctions $\set{e_k}$ and the fact that the Wiener process $W$ is real. Finally, the $\ob{\Delta x}^{-1}$ scaling in the noise term of (\ref{e42}) can be understood by taking the discrete Fourier transform (\ref{e15a}) of the noise term in discrete direct space, as introduced in (\ref{e17}).

\subsection{Decoupled Stochastic Allen-Cahn Equation}\label{numsimpli2}
The spatially decoupled stochastic Allen-Cahn equation (\ref{e16}) is conveniently discretized in direct space. Introducing the double-well intensity $g>0$ and the notation $x_j=\frac{2\pi}{N} j$,  where $j=(j_1, j_2)$ and $j_i=-N/2,\ldots,N/2-1$, we obtain the discretized version
\begin{align}\label{e43}
du_N(x_i)=g\cb{u_N(x_j)-u_N^3(x_j)}dt + \sigma dW_N(x_j).
\end{align}
Again, an implicit-explicit time stepping scheme was employed for (\ref{e43}), combining the explicit Euler-Maryama step for nonlinearity and noise with the trapezoidal scheme for the linear contribution
\begin{align}\label{e44}
u_N^{m+1}(x_j)=u_N^m(x_j)+ g\,\frac{\Delta t}{2} \cb{ u_N^{m+1}(x_j)+u_N^{m}(x_j)} - g \Delta t \cb{u_N^m(x_j)}^3 + \sigma \frac{\sqrt{\Delta t}}{\Delta x}\xi_j^m,
\end{align}
where $\{\xi_j^m\}$ are real i.i.d.\ $\mathcal{N}(0,1)$ random variables.

\subsection{Stochastic Allen-Cahn equation}\label{numAC}
Since the stochastic Allen-Cahn equation
\begin{align}\label{e45}
du= \cb{ \alpha \Delta u + g\ob{u-u^3}} dt + \sigma dW
\end{align}
involves the characteristic features of both prior models, the diffusion term as well as the nonlinearity, there is no natural choice between direct and Fourier space. Instead, we formulate the scheme in Fourier space, go to direct space using the inverse fast Fourier transform (IFFT), process the nonlinearity $u_N\mapsto u_N^3$, and go back to Fourier space using the fast Fourier transform (FFT). More precisely, the discretization scheme for (\ref{e45}) reads
\begin{align}\label{e46}
\hat{u}^{m+1}_N(k)= \hat{u}^m_N(k) + (g-\alpha \abs{k}^2) \frac{\Delta t }{2} \cb{\hat{u}^{m+1}_N(k)+\hat{u}^m_N(k)} - g\,\Delta t \ip{\ob{u_N^m}^3, e_k}_N +  \sigma \frac{\sqrt{\Delta t}}{\Delta x}\xi_k^m,
\end{align}
where $\xi_k^m=\frac{1}{\sqrt{2}}\ob{\eta_k^{m}+ i\, \zeta_k^{m}}$ with $\set{\eta_k^{m}}$ and $\set{\zeta_k^m}$ real i.i.d.\ $\mathcal{N}(0,1)$ random variables. In (\ref{e46}) we have used the notation
\begin{align}\label{e47}
\ip{x^3,e_k}_N:= \cb{\text{FFT}\cb{\ob{\text{IFFT} \cb{x}}^3}}_k,
\end{align}
for $x\in \text{span}\set{e_k :  k_{1,2}=-N/2,\ldots,N/2-1}$. Straight-forward implementation of the procedure (\ref{e47}) gives raise to an aliasing error, and we employed the so-called {\it Two-Thirds Rule} \cite{patterson1971spectral} to correct for this error. We refer to Appendix A in \cite{giada2002pseudospectral} for a detailed discussion of aliasing issues and the {\it Two-Thirds Rule}.

\section{Conclusions and Outlook}\label{conclusions}

\noindent Let us now discuss the implications of our work with respect to the discrepancy of frequent use (applied sciences) and virtual neglect (mathematical community) of the white noise-driven Allen-Cahn equation in two dimensions. Whereas the well-posedness of the equation {\it per se} remains an open question, Conjecture \ref{con2} suggests that modeling with the stochastic Allen-Cahn equation is indeed problematic: if the mesh size in simulations was shrunk, the numerical solutions would converge to the zero-distribution.
 As illustrated in Fig.\ \ref{fig:FIG_4} there is no pattern formation in the continuum limit -- the zero-distribution is a trivial solution. Therefore, even if one were able to make sense of the continuum equation, the employed numerical schemes could not be used to approximate it. In addition to the pseudospectral method described in Section \ref{numAC}, we positively tested our conjecture for a finite difference method. Since the pathologies do not seem to be caused by subtle deficiencies in the numerical schemes, but rather by the roughness of white noise itself, we expect Conjecture \ref{con2} to hold for all standard numerical schemes. \\

\noindent  One way to interpret simulations of the stochastic Allen-Cahn equation is to view them as numerical approximations of equations driven not by white noise, but by a noise field having a finite correlation length $\lambda$.
Not only does this render the model equations well-posed, but it can also provide a justification of the numerical experiments: in fact, if $\lambda\ll\Delta x$, where $\Delta x$ is the smallest mesh size used in the simulations, the noise is indeed uncorrelated on the grid, and the white noise
simulations are reasonable. However, such a correction may be problematic as $\lambda$ is dictated by the physics of the problem: it might be difficult to determine the correct correlation length, or the latter might be too big ($\lambda>\Delta x$), in which case the white noise simulations are indeed invalid. \\

\noindent For future endeavors, our work leads to the following questions. First, Conjecture \ref{con2} poses a challenge in terms of conceiving a rigorous proof, and we will address this issue in a forthcoming publication \cite{hairer2011}. Second, since the irregularity of white noise increases with the spatial dimension, we expect equations in three and more dimensions to be ill-posed, too. But a generalization of Conjecture \ref{con2} for $d\geq 3$ is yet to be established. Finally, the stochastic Allen-Cahn equation is only one specific example of the class of white noise-driven nonlinear SPDEs that are used to model physical phenomena in two and three space dimensions. It remains to be established whether other models (see Introduction for a detailed list) suffer from similar pathologies.

\section{Acknowledgments}

MDR is grateful to M. Hairer and H. Weber for fruitful discussions regarding heuristic arguments; in particular, M. Hairer conceived the idea of rewriting (\ref{e27}) as (\ref{e28}) and employing the Wick product. NN and PFT hold Canada Research Chairs in Applied Mathematics. This study was supported by the Natural Sciences and Engineering Research Council of Canada.

\section{Appendix}\label{appA}
To prove Theorem \ref{prop3} we need the following lemma.
\begin{lemma}\label{lemma1} Let $V(x) = \frac{1}{4}x^4-\frac{1}{2}x^2$ be the double-well potential, and let $d\geq 1$. Then the stationary solution $u_N$ of the regularized problem (20)  has the following properties:
\begin{enumerate}[(i)]
\item $\mathbb{E}\, u_N(x_j)=0,\quad\forall j_{1,2}=-N/2, \ldots, N/2-1, \quad \forall N\geq 2$.
\item $\mathbb{E}\,|u_N(x_j)|^2 \sim  \sigma N^{d/2}, \quad \forall j_{1,2}=-N/2, \ldots, N/2-1, \quad \forall N\geq 2$.
\end{enumerate}
 \end{lemma}

\begin{proof}
(i) follows from the symmetry of the potential $V$. For (ii), we recall that the stationary solution of (\ref{e17}) has a distribution at every grid point given by 
\begin{equation}  p_N(x) = \frac{1}{\mathcal{N}} \exp\left\{ -\frac{2V(x)}{\sigma^2\rho^{d/2}} \right\} = \frac{1}{\mathcal{N}} \exp\left\{ -\frac{(x^4/2-x^2)}{\sigma^2\rho^{d/2}} \right\}, \qquad \rho = \frac{N^2}{4\pi^2}. \end{equation}
To simplify calculations, we denote
\begin{equation}\label{constant} c:=\frac{1}{2\sigma^2 \rho^{d/2}}=\frac{(2\pi)^d}{2\sigma^2} N^{-d}.\end{equation}Then $$ p_N(x) =\frac{1}{\mathcal{N}}\, e^{c}\,\exp\left\{ -c(x^2-1)^2 \right\}
$$
\begin{align}
\mathbb{E}|u_N(x_i)|^2&=\frac{\int_{\mathbb{R}} x^2 p_N(x) \,dx}{\int_{\mathbb{R}} p_N(x) \,dx}= \frac{\int_0^\infty x^2 p_N(x) \,dx}{\int_0^\infty  p_N(x) \,dx}  = \frac{\int_0^\infty x^2 \exp\left\{ -c(x^2-1)^2 \right\}\,dx }{\int_0^\infty \exp\left\{ -c(x^2-1)^2 \right\} \,dx } =:R(c),\end{align}  since the integrands are even. 
With the change of variables $z=\sqrt{c}x^2$ we get
  \begin{equation}\label{ratio} R(c) = \frac{\int_0^\infty \frac{\sqrt{z}}{\sqrt{c}} \exp\{-(z-\sqrt{c})^2\} \,\frac{1}{2c^{1/4}} dz}{ \int_0^\infty \exp\{-(z-\sqrt{c})^2\} \,\frac{1}{2\sqrt{z}c^{1/4}} dz}. 
  \end{equation}
Now we define 
  \begin{align} P(c):= \frac{\int_0^\infty z^{1/2} \exp\{ -(z-\sqrt{c})^2\} \, dz}{\int_0^\infty z^{-1/2} \exp\{ -(z-\sqrt{c})^2\}  dz} \end{align}
and use \eqref{constant} and \eqref{ratio} to see that
\begin{equation}\label{forfinal}R(c) =  \sqrt{2}\sigma \rho^{d/4} P(c).\end{equation}
Next we establish that $P(c)$ is bounded on $ c\leq \frac{\ob{2\pi}^2}{2\sigma} \iff N\geq 1$. We recall the integral representation of the modified Bessel function $K_\nu(z)$ (\cite{watson1995treatise}, p183)
$$ K_\nu(z)=\frac{1}{2}(\frac{z}{2})^{\nu}\int_0^\infty\, \frac{1}{t^{\nu+1}}\,\,\exp{\ob{ -t-\frac{z^2}{4t}}}  \,dt,$$
and use it together with the software {\tt Maple} to get 
\begin{equation}\label{Pc}
P(c) = \frac{\sqrt{c}}{2 K_{1/4}(\frac{c}{2})} \left( K_{3/4}(\frac{c}{2}) - K_{1/4}(\frac{c}{2}) \right). \end{equation}
Exploiting recursive differentiation formulae for Bessel functions \cite{watson1995treatise}, one can deduce from (\ref{Pc}) that $P(c)$ is bounded and decreasing on intervals $[0,M]$, for all $M>0$. Recalling \eqref{constant} we conclude that $P$ as a function of $N$ is increasing and uniformly bounded for $N\geq1$. Together with (\ref{forfinal}) this yields (ii).  
\end{proof}

\bibliography{JCP_Ryser}
	\bibliographystyle{plain}

\end{document}